\documentclass[11pt]{article}

\title{\sc\Large Conjugacy, orbit equivalence and classification of measure preserving\\ group actions}
\author{Asger T\" ornquist}

\usepackage{amsmath}
\usepackage{amstext}
\usepackage{amssymb}
\usepackage{amsthm}
\usepackage{enumerate}
\usepackage{mathrsfs}

 \DeclareMathOperator{\Aut}{Aut}

\DeclareMathOperator{\Irr}{Irr}\DeclareMathOperator{\Mod}{Mod}

\DeclareMathOperator{\Ind}{Ind}

\newcommand{\erg}{\mathcal E}
\newcommand{\wmix}{\mathcal W}
\newcommand{\act}{\mathcal A}

\DeclareMathOperator{\tfa}{TFA}

\def\R{{\mathbb R}}
\def\C{{\mathbb C}}

\def\N{{\mathbb N}}
\def\Z{{\mathbb Z}}

\def\({{\normalfont (}}
\def\){{\normalfont )}}

\newcounter{mypar}
\newcounter{thmcounter}

\newcommand{\mysec}[1]{\setcounter{thmcounter}{0}\addtocounter{mypar}{1}\section*{\begin{center}\normalsize{\sc \S \arabic{mypar}.
#1}\end{center}}}

\newcommand{\mythm}[2]{\addtocounter{thmcounter}{1}\subparagraph{{\sc \arabic{mypar}.\arabic{thmcounter}. #1}}{\it
#2}}

\newcommand{\thm}[2]{{\sc#1\ }{\it
#2}}

\newcommand{\defn}{\addtocounter{thmcounter}{1}\subparagraph{{\rm \arabic{mypar}.\arabic{thmcounter}.}
{\it Definition.}}}

\newcommand{\myheadpar}[1]{\addtocounter{thmcounter}{1}\subparagraph{\rm \arabic{mypar}.\arabic{thmcounter}. {\it #1}}}

\newcommand{\mypar}{\addtocounter{thmcounter}{1}\subparagraph{\rm \arabic{mypar}.\arabic{thmcounter}.}}

\newcommand{\claim}[1]{\subparagraph{{\sc #1}}}

\begin{document}

\maketitle

\begin{abstract}
We prove that if $G$ is a countable discrete group with property (T)
over an infinite subgroup $H\leqslant G$ which contains an infinite
Abelian subgroup or is normal, then $G$ has continuum many orbit
inequivalent measure preserving a.e. free ergodic actions on a
standard Borel probability space. Further, we obtain that the
measure preserving a.e. free ergodic actions of such a $G$ cannot be
classified up to orbit equivalence be a reasonable assignment of
countable structures as complete invariants. We also obtain a
strengthening and a new proof of a non-classification result of
Foreman and Weiss for conjugacy of measure preserving ergodic, a.e.
free actions of discrete countable groups.
\end{abstract}

\mysec{Introduction}

\paragraph{(A)} Let $G_1$ and $G_2$ be discrete countable groups, acting by measure
preserving transformations on standard Borel probability spaces
$(X,\mu)$ and $(Y,\nu)$, respectively, and giving rise to orbit
equivalence relations $E_{G_1}$ and $E_{G_2}$. Recall that the
actions of $G_1$ and $G_2$ are {\it orbit equivalent} if there is a
measure preserving transformation $T:X\to Y$ such that
$$
xE_{G_1} y\iff T(x) E_{G_2} T(y).
$$
The notion of orbit equivalence was introduced by Dye in
\cite{dye},\cite{dye2}, who proved that any two ergodic $\Z$-actions
are orbit equivalent. Connes, Feldman and Weiss generalized this in
\cite{cfw} to all countable amenable groups, and eventually the
combined work of Connes, Weiss \cite{connesweiss}, Schmidt
\cite{schmidt}, and Hjorth \cite{hjorth2} showed that this
characterizes amenability: Any non-amenable discrete group has at
least two non-orbit equivalent ergodic a.e. free measure preserving
actions on a standard Borel probability space.

\medskip

Recall that a group $G$ is said to have the {\it relative property }
(T) over a subgroup $H\leqslant G$ if there is $Q\subseteq G$ finite
and $\varepsilon>0$ such that whenever $\pi:G\to\mathscr X$ is a
unitary representation of $G$ on a Hilbert space $\mathscr X$ with a
unit vector $\xi\in\mathscr X$ such that for all $g\in Q$
$$
\|\pi(g)(\xi)-\xi\|_{\mathscr X}<\varepsilon,
$$
then there is a non-zero $H$-invariant vector in $\pi$. The purpose
of this paper is to prove the following (Theorem 4, \S4):

\medskip

\mythm{Theorem.}{If $G$ is a countable discrete group with the
relative property (T) over an infinite subgroup $H\leqslant G$ such
that either

\(i\) $H$ contains an infinite Abelian subgroup \emph{or}

\(ii\) $H$ is normal in $G$,

\noindent then $G$ has continuum many orbit inequivalent a.e. free
measure preserving actions on a standard Borel probability space.}

\bigskip

Such groups include well known examples such as the pair
$(SL_2(\Z)\ltimes\Z^2,\Z^2)$, as well as any group with Property (T)
outright, and groups that are formed by operations that preserve the
classes (i) and (ii), such as free and direct products of arbitrary
countable groups with a group with the relative property (T) over an
infinite abelian subgroup.

It should be noted that (ii) in Theorem 1.1 was proven by Popa in
\cite{popa}. We obtain a different proof here. The condition of
normality can be weakened, see \S3 Theorem 3.

\paragraph{(B)} In descriptive set theory, a theory of
classification has emerged over the last two decades, centered
around the notion of {\it Borel reducibility}: If $X$ and $Y$ are
Polish spaces and $E$ and $F$ are equivalence relations on $X$ and
$Y$, respectively, we say that $E$ is Borel reducible to $F$ if
there is a Borel $f:X\to Y$ such that
$$
(\forall x,y\in X) xEy\iff f(x) F f(y).
$$
We write $E\leq_B F$ if this is the case, and $E<_B F$ if $E\leq_B
F$ but not $F\leq_B E$. Considering $E$ and $F$ as (potentially
interesting) classification problems among the points of $X$ and
$Y$, respectively, we can interpret the statement $E\leq_B F$ as
saying that the points of $X$ are classified up to $E$ equivalence
by a Borel assignment of $F$-equivalence classes as complete
invariants. Alternatively, we may think of $E\leq_B F$ as an
expression of the relative difficulty posed by the two
classification problems. The reason that $f$ is required to be Borel
is to guarantee that the classification is reasonably
``computable'', thus avoiding the useless classifications offered by
the axiom of choice.

The notion of Borel reducibility will allow us to state the theorems
of this paper in a much sharper way. We say that an equivalence
relation $E$ (on a Polish space $X$) is {\it smooth} or {\it
concretely classifiable} if there is a Borel $f:X\to\R$ such that
$$
(\forall x,y\in X) xEy\iff f(x)=f(y),
$$
that is, if there is a complete classification using real numbers as
invariants. A pivotal example of a {\it non-smooth} equivalence
relation is $E_0$, defined on $2^\N$ by
$$
x E_0 y\iff (\exists N)(\forall n\geq N) x(n)=y(n).
$$
However, there are natural equivalence relations far beyond the
complexity of $E_0$. An important divide in the landscape of
classification is the distinction between those equivalence
relations that allow classification by {\it countable structures},
that is, allow a Borel assignment of invariants such as countable
groups, graphs, trees, etc., considered up to the natural notion of
isomorphism, as oppose to those that do not. (Precise definitions
are given in \S 2 below.) Then we have the following sharper version
of Theorem 1.1:

\medskip

\mythm{Theorem.}{If $G$ is a countable discrete group with the
relative property (T) over an infinite subgroup $H\leqslant G$ such
that either

\(i\) $H$ contains an infinite Abelian subgroup \emph{or}

\(ii\) $H$ is normal in $G$,

Then orbit equivalence of measure preserving ergodic $G$ actions is
not smooth and cannot be classified by countable structures.}

\bigskip

It should be noted that by a result in \cite{atphd}, which was based
on the result of Popa in \cite{popa}, that in the case (ii) above it
was known that the isomorphism relation for countable torsion free
abelian groups $\simeq_{\tfa_{\aleph_0}}$ is Borel reducible to
orbit equivalence. Thus the above Theorem shows that it is in fact
strictly above $\simeq_{\tfa_{\aleph_0}}$ in the $\leq_B$-hierarchy,
and rules out that countable discrete groups can, in a reasonable
way, be used as complete invariant for orbit equivalence of groups
with the relative property (T).

We refer the reader to the introduction of \cite{hjorth3} for
further background on the very interesting subject of Borel
reducibility and the descriptive set theory of classification
problems.

\paragraph{(C)} The proof of Theorem 1.2 (and 1.1)
is based on the method used by Hjorth in \cite{hjorth2}, relating
{\it conjugacy} of measure preserving actions to orbit equivalence.
The most important case is when $G$ is a group with the relative
property (T) over an infinite Abelian subgroup $H\leqslant G$.

Let $\sigma$ and $\tau$ be measure preserving $G$ actions on a
standard Borel probability space $(X,\mu)$. We will say that
$\sigma$ and $\tau$ are \emph{conjugate on $H$} if $\sigma|H$ and
$\tau|H$ are conjugate. We write $\sigma\simeq_H\tau$. We say that
$\sigma$ is \emph{ergodic on $H$} if $\sigma|H$ is ergodic. In \S 4
we prove the following relativization of Corollary 2.6 in
\cite{hjorth2}:

\mythm{Lemma.}{Let $G$ be a countable group with the relative
property (T) over the infinite subgroup $H\leqslant G$. Then at most
countably many ergodic-on-$H$ probability measure preserving actions
of $G$ that are not conjugate on $H$ can be orbit equivalent.}

\bigskip

Given this, the task becomes to produce many $G$ actions that are
non-conjugate and ergodic on $H$. To this end we prove a version of
Foreman and Weiss anti-classification result  in \cite{foremanweiss}
(Theorem A) for extensions of actions. This is based on Theorem 5.1
of Kechris and Sofronidis in \cite{kecsof}:

\bigskip

\thm{Theorem ({\rm Kechris-Sofronidis}).}{Let $X$ be an uncountable
Polish space and let $P_c(X)$ be the Polish space of continuous
Borel probability measures. Then absolute equivalence $\approx$ in
$P_c(X)$ is generically turbulent\footnote{See \S2 for a
definition.} and thus it is not smooth and not classifiable by
countable structures.}

\bigskip

Assuming now that $H\leqslant G$ is an infinite Abelian group, then
for each symmetric measure $\sigma$ on the dual $\hat H$ we have an
associated real positive definite function
$$
\varphi_\sigma(h)=\int \chi(h)d\sigma(\chi).
$$
Then if $\bar\varphi_\sigma$ is the trivial extension to $G$, then
we may construct a Gaussian action $\rho_\sigma$ of $G$ as in
\cite{glasner} p. 90-91. In Lemma 3.5 we will compute a formula for
the maximal spectral type of the Koopman representation
$\kappa^\sigma|H$ of this action restricted to $H$. We then show
that if the measure $\sigma$ has support on a ``sufficiently
independent'' set (see Definition 3.11), then
$\sigma\mapsto\rho_\sigma$ provides us with a Borel function such
that equivalent measures produce conjugate $G$-actions, and
inequivalent measures produce $G$-actions that are non-conjugate on
$H$. Furthermore, the action $\rho_\sigma|H$ will be weakly mixing.
Thus $\sigma\mapsto\rho_\sigma$ is a Borel reduction of $\approx$ to
$\simeq_H$.

\medskip

We note that we obtain a new proof and a stronger version of Foreman
and Weiss' result (Theorem 1, 2 and 3, \S 3):

\mythm{Theorem.}{Let $G$ be a countable discrete group and let
$H\leqslant G$ be an infinite subgroup. Let $(X,\mu)$ be a standard
Borel probability space. Then:

\(i\) The measure preserving $H$-actions on $(X,\mu)$ that can be
extended to an ergodic measure preserving $G$-action on $(X,\mu)$
cannot be classified up to conjugacy by countable structures.

\(ii\) If $H$ either contains an infinite Abelian subgroup or is
normal in $G$, then the measure preserving ergodic $H$-actions on
$(X,\mu)$ that can be extended to an ergodic measure preserving
$G$-action on $(X,\mu)$ are not classifiable by countable
structures.}

\paragraph{(D)} The author would like to thank Greg Hjorth for
useful discussions about the ideas for this paper at an early stage
while the author was visiting the Mathematics Department and Logic
Center at UCLA, in July 2006, and Alexander Kechris for making
available his extensive new notes on measure preserving actions
\cite{kechris2}, the appendices of which were indispensable to this
work. The author would also like to thank the referee, whose
comments and suggestions have improved this paper vastly.

Research for this paper was supported in part by the Danish Natural
Science Research Council grant no. 272-06-0211.

\mysec{Borel reducibility}

\myheadpar{Notation:} Throughout this paper, if $E$ is an
equivalence relation on some set $X$ and $Y\subseteq X$, we will
write $E^Y$ for $E|Y\times Y$.

\myheadpar{Reducibility notions.} We first introduce two variant
technical notions of Borel reducibility that will be of useful for
us in this paper. Let $X,Y$ be Polish spaces, and let $E$ and $F$ be
equivalence relations on $X$ and $Y$, respectively. Recall that we
say that $E$ is Borel reducible to $F$, written $E\leq_B F$, if
there is a Borel $f:X\to Y$ such that
$$
(\forall x,y\in X) x E y\iff f(x)Ff(y).
$$
If $F'$ is another equivalence relation on $Y$ and $F\subseteq F'$
then it is natural to introduce the following stronger notion: We
will say that $E$ is Borel reducible to the pair $(F,F')$, written
$E\leq_B (F,F')$, if there is a Borel function $f:X\to Y$ such that
$$
(\forall x,y\in X) xEy\implies f(x)Ff(y)\wedge f(x)F'f(y)\implies
xEy.
$$
It is clear that if $E\leq_B (F,F')$ then for any equivalence
relation $F\subseteq F''\subseteq F'$ we have $E\leq_B F''$.

A useful weaker reducibility notion is the following: We say that
$E$ is {\it Borel countable-to-one reducible} to $F$, written
$E\leq_B^\omega F$, if there is a Borel $f:X\to Y$ such that
$$
xEy\implies f(x) F f(y)
$$
and the induced map $\hat f: X/E\to Y/F:\hat f([x]_E)=[f(x)]_F$ is
countable-to-1.

\myheadpar{Classification by countable structures.} Recall that an
equivalence relation $E$ {\it admits a classification by countable
structures} (or {\it is classifiable by countable structures}) if
there is a countable language $\mathcal L$ such that
$E\leq_B\simeq_{\Mod(\mathcal L)}$, where $\simeq_{\Mod(\mathcal
L)}$ is isomorphism for countable models of $\mathcal L$ (see
\cite{hjorth3} for details.) Let $E$ be an equivalence relation on a
Polish space $X$. We will say that $E$ is {\it generically
turbulent} if
\begin{enumerate}[i)]

\item All $E$ classes are meagre and

\item There is a Polish group $G$ acting continuously on $X$ such
that the induced orbit equivalence relation $E_G\subseteq E$ and the
action of $G$ on $X$ is generically turbulent in the sense of
\cite{hjorth3}.
\end{enumerate}

Then it is immediate from \cite{hjorth3}, Theorem 3.18, and
\cite{beke}, Theorem 3.4.5, that if $E$ is generically turbulent
then $E_0\leq_B E$ and $E$ is not classifiable by countable
structures. Further we have:

\mythm{Proposition.}{Let $F$ be an equivalence relation on a Polish
space $Y$, and suppose there is a generically turbulent equivalence
relation $E$ on a Polish space $X$ such that $E\leq_B^\omega F$.
Then $E_0<_B F$ and $F$ does not admit a classification by countable
structures.}

\begin{proof}
Let $f$ witness that $E\leq_B^\omega F$. Define
$$
xE' y\iff f(x)Ff(y).
$$
Then $E\subseteq E'$ and $E'$ has countable index over $E$. Thus
$E'$ is generically turbulent. Since by definition $E'\leq_B F$ it
follows that $E_0\leq_B F$ and that $F$ is not classifiable by
countable structures.
\end{proof}

\mysec{Conjugacy on a subgroup}

Let $(X,\mu)$ be a standard Borel probability space. We denote by
$\Aut(X,\mu)$ the group of measure preserving transformations on
$X$, which we give the weak topology (c.f. \cite{kechris}). Further,
if $G$ is a countable group, we denote by $\act(G,X,\mu)$ the set of
homomorphisms of $G$ to $\Aut(X,\mu)$ equipped with the subspace
topology it inherits from $\Aut(X,\mu)^G$, and identify
$\act(G,X,\mu)$ with the space of measure preserving
transformations. $\erg(G,X,\mu)$ and $\wmix(G,X,\mu)$ denote the
subsets of ergodic, respectively weakly mixing, measure preserving
actions of $G$ on $X$. We denote by $\act^*(G,X,\mu)\subseteq
\act(G,X,\mu)$ subset of a.e. free actions, and we let
\begin{align*}
&\erg^*(G,X,\mu)=\erg(G,X,\mu)\cap \act^*(G,X,\mu)\\
&\wmix^*(G,X,\mu)=\wmix(G,X,\mu)\cap \act^*(G,X,\mu).
\end{align*}
$\Aut(X,\mu)$ acts on $\act(G,X,\mu)$ by conjugation. We denote by
$\simeq$ the conjugation relation in $\act(G,X,\mu)$.

If $H$ is a subgroup of $G$, then we denote by $\erg_H(G,X,\mu)$ the
subset of $\sigma\in\erg(G,X,\mu)$ such that $\sigma|H$ is ergodic
as an action of $H$, that is, the ergodic on $H$ m.p. actions of $G$
as defined in the introduction. Similarly, we denote by
$\wmix_H(G,X,\mu)$ the set of $\sigma\in\wmix(G,X,\mu)$ such that
$\sigma|H$ is weakly mixing. We denote by $\simeq_H$ the equivalence
relation in $\act(G,X,\mu)$ of being conjugate on $H$: That is, for
$\sigma,\tau\in\act(G,X,\mu)$ we write $\sigma\simeq_H\tau$ iff
there is $T\in\Aut(X,\mu)$ such that on a measure one set we have
$$
(\forall h\in H) T(\sigma(h)(x))=\tau(g)(T(x)).
$$
Finally, we define the corresponding sets of free actions
$\erg_H^*(G,X,\mu)$ and $\wmix_H^*(G,X,\mu)$ as before.

Let $X$ be a Polish space. We will denote by $P(X)$ the Polish space
of all Borel probability measures, c.f. \cite{kechris} pp. 109. We
denote by $P_c(X)$ the subset of continuous (i.e. non atomic)
probability measures. If $X$ is uncountable then $P_c(X)$ is a dense
$G_\delta$ subset of $P(X)$, cf. \cite{kecsof}. Absolute equivalence
$\approx$ is a Borel equivalence relation in $P(X)$, c.f.
\cite{kechris} 17.39. In this section we will prove the following:

\medskip

\thm{Theorem 1.}{Let $G$ be a countable group and $A\leqslant G$ an
infinite abelian subgroup, and let $(X,\mu)$ be a non-atomic
standard Borel probability space. Then
$\approx^{P_c(2^\N)}\leq_B(\simeq^{\wmix^*_{A}(G,X,\mu)},\simeq_A^{\wmix^*_{A}(G,X,\mu)})$.
So for any $A\leqslant H\leqslant G$ we have that
$E_0<_B\simeq_H^{\wmix^*_{H}(G,X,\mu)}$ and
$\simeq_H^{\wmix^*_{H}(G,X,\mu)}$ is not classifiable by countable
structures.}

\bigskip

The proof is presented in a sequence of lemmata. For the rest of
this section, $A\leqslant G$ is a fixed infinite abelian subgroup of
$G$, and $G$ is a fixed countable group. We also fix once and for
all a sequence $(g_i)$ of coset representatives, and suppose that
the identity is chosen to represent $A$.

\medskip

The idea of the proof is the following: Given a finite (continuous)
measure $\sigma$ on the dual $\hat A$, we have the corresponding
unitary $\pi^\sigma$ on $L^2(\hat A,\sigma)$. If $\sigma$ is a
symmetric measure, then the corresponding positive definite function
$\varphi_\sigma$ is real valued. Considering then
$\Ind_A^G(\pi^\sigma)$ of $G$, which corresponds to the trivial
extension $\bar\varphi_\sigma$ of $\varphi_\sigma$, we can then
obtain a Gaussian $G$-action $\rho_\sigma$ as in \cite{glasner} (p.
90-91). Considering this action as an $A$-action, we will compute
the maximal spectral type of the corresponding Koopman
representation, using a notion of ``lifted'' spectral measure,
introduced below. We will then show that if the measures
$\sigma,\sigma'$ have support on a suitably chosen perfect subset
$C$ of $\hat A$, then the resulting Gaussian action $\rho_\sigma|A$
and $\rho_{\sigma'}|A$ are conjugate as $A$-actions if and only if
$\sigma$ and $\sigma'$ are absolutely equivalent.
\medskip

\myheadpar{Lifts.} Let $H<A$ be a subgroup of the abelian group $A$.
We define $P_H:\hat A\to\hat H$ by
$$
P_H(\chi)=\chi|H.
$$
Let $\sigma$ be a positive finite measure on $\hat H$, and let
$\tau:\hat H\to P(\hat A):\chi\mapsto\tau_\chi$ be a Borel
assignment of probability measures such that for all $\chi\in\hat H$
we have
$$
\tau_\chi(P_H^{-1}(\chi))=1.
$$
Then we define the {\it $\tau$-lift} of $\sigma$ to be the measure
$\sigma_{\tau}$ on $\hat A$ defined by
$$
\int f d\sigma_{\tau}=\iint f(\eta) d\tau_\chi(\eta)d\sigma(\chi).
$$
Note that the lifting process respects absolute continuity: if
$\sigma\ll\sigma'$ then $\sigma_\tau\ll\sigma'_\tau$. Note also that
if $H'\leqslant H$ is a subgroup, then
$P_{H'}[\sigma_\tau]=P_{H'}[\sigma]$. In particular,
$P_{H}[\sigma_\tau]=\sigma$.

\bigskip

\myheadpar{Haar lifts.} A special instance of this construction
arises as follows: Identify $\widehat{A/H}$ with the set of
characters $\chi\in \hat A$ such that $\chi|H=1$, and let
$\lambda_{\widehat {A/H}}$ be the (normalized) Haar measure on
$\widehat{A/H}$. Then for $\chi\in \hat H$ define for a Borel
$B\subseteq \hat A$,
$$
\lambda_\chi(B)=\lambda_{\widehat{A/H}}({\tilde\chi}^*B),
$$
where $\tilde\chi$ is some (any) extension of $\chi$ to a character
on $\hat A$.\footnote{We note that we may find a Borel
$\chi\mapsto\tilde\chi$ function for this purpose since we may find
a transversal for the equivalence relation induced by the natural
action of the compact group $\widehat{A/H}$ on $\hat A$, using
\cite{kechris2}, Theorem 12.16.} Note that $\lambda_\chi$ is
independent of the choice of extension $\tilde\chi$ since
$\lambda_{\widehat{A/H}}$ is translation invariant. Then
$\lambda_\chi(P_H^{-1}(\chi))=1$ and $\chi\mapsto\lambda_\chi$ is a
Borel assignment of Borel probability measures on $\hat A$. For a
measure $\sigma$ on $\hat H$ we call the lifted measure
$\sigma_{\lambda}$ on $\hat A$ the {\it Haar lift} of $\sigma$, and
denote it $\sigma_H^A$, with a purposely suggestive notation. It is
routine to check that if $\varphi_\sigma:H\to \C$ is the positive
definite function corresponding to the spectral measure $\sigma$ on
$\hat H$, then $\sigma^A_H$ is the spectral measure of the trivial
extension of $\varphi_\sigma$ to $A$.

\bigskip

For each coset representative $g_i$, define $H_i=A\cap
g_{i}Ag_{i}^{-1}$ and $H'_i=g_{i}^{-1}Ag_{i}\cap A$. We also define
$\varphi_i:H_i\to H'_i$ to be the isomorphism
$$
\varphi_i(a)=g_i^{-1}ag_i.
$$
Note that $\varphi_i$ only depends on the coset, not on the choice
of representative. We denote by $\hat\varphi_i:\hat H_i'\to \hat
H_i$ the associated isomorphism of the character groups,
$\hat\varphi_i(\chi)=\chi\circ\varphi_i$. More generally, and in a
slight abuse of notation, we will also denote by $\hat\varphi_i$ the
homomorphism $\hat A\to\hat H_i$ defined by
$\hat\varphi_i(\chi)=\chi\circ\varphi_i$. If there is any danger of
confusion, that latter will be written $\hat\varphi_i\circ
P_{H_i'}$. We leave the proof of the next Lemma as an easy exercise
for the reader in using the definition of lifts and induced
representation:

\mythm{Lemma.}{Suppose $(\pi,\mathscr X)$ is a unitary
representation of $A$ on a separable Hilbert space $\mathscr X$ and
that $A$ is a subgroup of some countable discrete group $G$. If the
maximal spectral type of $\pi$ is $\sigma$, then the maximal
spectral type of $\Ind_A^G(\pi)|A$ is
$$
\sigma_{\Ind_A^G(\pi)|A}^{\max}\approx \sum_i \frac 1 {2^i}
\hat\varphi_i[\sigma]^A_{H_i}.
$$
}

\mypar Let $\sigma$ be a symmetric, finite and positive measure on
$\hat A$, so that the associated positive definite function
$\varphi_\sigma$ is real valued. Let $\bar\varphi_\sigma$ be the
trivial extension of $\varphi_\sigma$ to $G$, and let
$\mu_\sigma=\mu_{\bar\varphi_\sigma}$ be the associated Gaussian
measure on $\R^G$, which we denote $\rho_\sigma$. If $\kappa$ is the
associated Koopman representation on $L^2(\R^G,\mu)$ then
$$
\kappa^{H_1^{(c)}(\mu)}\sim_u \Ind_A^G(\pi^\sigma),
$$
where $H_1^{(c)}$ is the 1st Wiener Chaos. It follows from Lemma 3.3
that the maximal spectral type of $\kappa^{H_1^{(c)}(\mu)}|A$, which
we denote by $w_{\sigma}$, is
$$
w_\sigma=\sum_i \frac 1 {2^i} \hat\varphi_i[\sigma]^A_{H_i}.
$$
Thus the maximal spectral type of the Koopman representation
$\kappa$ is $\exp(w_\sigma)$, from which we obtain:

\medskip

\mythm{Lemma.}{With notation as in 3.4, the maximal spectral type of
$\kappa_0|A$ is
$$
\sum_{\substack{i_1,\ldots i_n\in\N\\ n\in\N}} \frac 1 {2^n
2^{i_1+\cdots +i_n}}
\hat\varphi_{i_1}[\sigma]^A_{H_{i_1}}*\cdots*\hat\varphi_{i_n}[\sigma]^A_{H_{i_n}}.
$$
}

\defn For a finite positive measure $\sigma$ on $\hat
A$, define
$$
\Phi_\sigma= \sum_{\substack{i_1,\ldots i_n\in\N\\ n\in\N}} \frac 1
{2^n 2^{i_1+\cdots +i_n}}
\hat\varphi_{i_1}[\sigma]^A_{H_{i_1}}*\cdots*\hat\varphi_{i_n}[\sigma]^A_{H_{i_n}},
$$
so that the previous Lemma states that the maximal spectral type of
$\kappa_0|A$ on $L^2_0(\R^G,\mu_\sigma)$ is $\Phi_\sigma$.

\bigskip

The following Lemma is routine and we omit the proof.

\mythm{Lemma.}{Let $H<A$ be a subgroup of $A$ with infinite index
and suppose $C\subseteq\hat A$ is a Borel set such that for $P_H|C$
is 1-1. Then:

\(i\) For any finite positive Borel measure $\sigma$ on $\hat H$ and
$\chi\in \hat A$ it holds for the Haar lift that $\sigma_H^A(\chi
C)=0$.

\(ii\) For any finite positive Borel measure $\sigma$ on $\hat H$
and any finite positive Borel measure $\tau$ on $\hat A$ we have
that
$$
\tau*\sigma_H^A(C)=0.
$$
}

\defn A subset $C\subseteq \hat A$ is called {\it sufficiently
independent} if

\begin{enumerate}

\item For all $i$ such that $H_i$ is infinite we have that $P_{H_i}|C$ is
1-1, and for all $i$ such that $H_i'$ is infinite we have that
$P_{H_i'}|C$ is 1-1.

\item For all $\chi_1,\ldots, \chi_k,\chi_{k+1}\in C$ distinct
and all $H_{i_1},\ldots H_{i_{n_k}}$ with $[A: H_{i_j}]<\infty$ for
all $j=1,\ldots, k$, and all $s:\{1,\ldots, k\}\to \{-1,1\}$ we have
that
$$
\prod_{j=1}^{k}{\chi_j^{s(j)}\circ\varphi_{i_j}}|H_{(i_1,\ldots,i_k)}\neq\chi_{k+1}|H_{(i_1,\ldots,i_k)},
$$
where
$$
H_{(i_1,\ldots,i_k)}=\bigcap_{j=1}^{k} H_{i_j}.
$$
\end{enumerate}

Moreover, for a measure $\sigma$ on $\hat A$ we define the {\it
symmetrization} $\hat \sigma$ of $\sigma$ by
$$
\hat\sigma(B)=\frac 1 2(\sigma(B)+\sigma(B^*)),
$$
where $B\subseteq\hat A$ is Borel. Note that if $\psi:\hat A\to H$
is a Borel homomorphism of $\hat A$ into some (Polish) group $H$
then $\psi[\hat\sigma]=\widehat{\psi[\sigma]}$.

\bigskip

The following is immediate from the definition of being sufficiently
independent:

\mythm{Lemma.}{Suppose $C$ is a sufficiently independent subset of
$\hat A$. If $\sigma$ is a continuous Borel measure on $C$ then for
any $i$ such that $H_i$ is infinite $P_{H_i}[\sigma]$ and
$P_{H_i'}[\sigma]$ are continuous. Moreover, the measure
$\hat\varphi_i[\sigma]$ is continuous.}

\bigskip

We then have:

\mythm{Lemma.}{Suppose $C$ is a perfect sufficiently independent
subset of $\hat A$ such that $\lambda_{\hat A}(C)=0$, and let
$\sigma,\tau$ be finite non-atomic positive Borel measures on $C$.
Then $\sigma\not\approx\tau$ implies
$\Phi_{\hat\sigma}\not\approx\Phi_{\hat\tau}$. }
\begin{proof}
Suppose $B\subseteq C$ is a Borel set such that $\sigma(B)>0$ and
$\tau(B)=0$. Note that by definition $\sigma\ll\Phi_{\hat\sigma}$,
thus $\Phi_{\hat\sigma}(B)>0$. For simplicity, we will write
$\sigma_i$ for $\hat\varphi_i[\sigma]_{H_i}^A$ and $\tau_i$ for
$\hat\varphi_i[\tau]_{H_i}^A$.

We claim that $\Phi_{\hat\tau}(B)=0$. For this, it suffices to show
that for all $i_1,\ldots i_n\in\N$ we have
$$
\hat\tau_{i_1}*\cdots*\hat \tau_{i_n}(B)=0.
$$
If for some $1\leq k\leq n$ we have that $[A:H_{i_k}]=\infty$ and
$|H_{i_k}|=\infty$ then it follows from Lemma 3.7 that
$$
\hat\tau_{i_1}*\cdots*\hat \tau_{i_n}(B)=0.
$$
Also, if $|H_{i_k}|<\infty$ for some $1\leq k\leq n$ is clear that
$\hat\tau_{i_k}\ll\lambda_{\hat A}$. It follows that\footnote{We use
the notation $\prod_{i=1}^n
*\mu_i=\mu_1*\cdots*\mu_n$.}
$$
\prod_{j=1}^n*\hat\tau_{i_j}\ll (\prod_{\substack{j=1\\ j\neq
k}}^n*\hat\tau_{i_j})*\lambda_{\hat A}=\lambda_{\hat A},
$$
where the latter follows from the translation invariance of the Haar
measure $\lambda_{\hat A}$. Since we assumed $\lambda_{\hat A}(C)=0$
it follows that
$$
\hat\tau_{i_1}*\cdots*\hat \tau_{i_n}(B)=0.
$$

So without loss of generality, assume that $[A:H_{i_j}]<\infty$ for
all $j=1,\ldots, n$. Let
$$
H=H_{(i_1,\ldots, i_k)}=\bigcap_{j=1}^n H_{i_j}.
$$
Then $[A:H]<\infty$, too. Thus $P_H$ is finite-to-1 and so $P_H(B)$
is Borel. We will show that
$$
P_H[\prod_{j=1}^n*\hat\tau_{i_j}](P_H(B))=0,
$$
from which it clearly follows that $\hat\tau_{i_1}*\cdots*\hat
\tau_{i_n}(B)=0$.

Since $P_H$ is a homomorphism we have
$$
P_H[\prod_{j=1}^n*\hat\tau_{i_j}]=\prod_{j=1}^n*
P_H[\hat\tau_{i_j}].
$$
Further, we have $P_H[\hat\tau_{i_j}]=P_H\circ\hat\varphi_{i_j}
[\hat\tau]$ since
$\hat\tau_{i_j}=\hat\varphi_{i_j}[\hat\tau]_{H_{i_j}}^A$.

Suppose for a contradiction that
$$
P_H[\prod_{j=1}^n*\hat\tau_{i_j}](P_H(B))=\prod_{j=1}^n*P_H\circ\hat\varphi_{i_j}[\hat\tau](P_H(B))>0.
$$
By definition then
$$
P_H\circ\hat\varphi_{i_1}[\hat\tau]\times\cdots\times
P_H\circ\hat\varphi_{i_n}[\hat\tau](\{(\chi_1,\ldots,\chi_n)\in\hat
H:\chi_1\cdots \chi_n\in P_H(B)\})>0.
$$
It follows by repeated applications of Fubini's Theorem we may then
find $\chi_1,\ldots,\chi_n\in C\setminus B$ distinct and
$s:\{1,\ldots, n\}\to \{-1,1\}$ such that
$$
P_H\circ \hat\varphi_{i_1}(\chi_1^{s(1)})\cdots P_H\circ
\hat\varphi_{i_n}(\chi_n^{s(1)})=P_H(\eta)
$$
for some $\eta\in B$. But this says precisely that
$$
\prod_{j=1}^{n}{\chi_j^{s(j)}\circ\varphi_{i_j}}|H=\eta|H,
$$
flatly contradicting that $C$ is sufficiently independent.
\end{proof}

\bigskip

\mythm{Lemma.}{There is a perfect sufficiently independent subset of
$\hat A$.}

\begin{proof}
We will use Mycielski-Kuratowski's Theorem (\cite{kechris}, Theorem
19.1, p. 129.) For each infinite subgroup $H<A$ let
$$
R_H=\{(\chi,\eta)\in\hat A\times\hat A: P_H(\chi)\neq P_H(\eta)\}
$$
and for each sequence $i_1,\ldots, i_n\in\N$ where  such that
$[A:H_{i_j}]<\infty$ for all $1\leq j\leq n$ and $ s:\{1,\ldots,
n\}\to\{-1,1\}$ let
$$
Q^{(i_1,\ldots, i_n)}_s=\{(\chi_1,\ldots,\chi_m,\eta)\in\hat
A^{m+1}:\prod_{j=1}^{n}\chi_j^{s(j)}\circ\varphi_{i_j}|H_{(i_1,\ldots,
i_n)}\neq\eta|H_{(i_1,\ldots, i_n)}\}
$$
where
$$
H_{(i_1,\ldots, i_n)}=\bigcap_{j=1}^n H_{i_j}.
$$

\claim{Claim 1.} For each infinite $H<A$ the relation $R_H$ is
comeagre in $\hat A\times\hat A$.

\bigskip

{\it Proof}. By Kuratowski-Ulam's Theorem it suffices to show that
for each $\eta\in\hat A$ the set
$$
\{\chi\in\hat A: P_H(\chi)=P_H(\eta)\}=\eta\{\chi\in\hat A:
\chi|H=1\}
$$
is meagre, and so to show that $\{\chi\in\hat A: \chi|H=1\}$ is
meagre. Suppose not. Then by Pettis' Theorem (\cite{kechris2}, 9.9)
it contains an open neighborhood of the identity. The map
$\chi\mapsto\chi|H$ is open by the open mapping Theorem, and so it
follows that $\{1\}$ is open in $\hat H$. Thus $\hat H$ must be
discrete, and so $H$ itself must be finite. \qed (Claim 1)

\claim{Claim 2.} $Q^{(i_1,\ldots, i_n)}_s$ is comeagre in $\hat
A^{n+1}$ .

\bigskip

{\it Proof}. Let $H=H_{(i_1,\ldots, i_n)}$. By Kuratowski-Ulam's
Theorem it is sufficient to show for each $\eta\in\hat A$ that the
set
$$
\{(\chi_1,\ldots,\chi_n)\in\hat
A^n:\prod_{j=1}^{n}\chi_j^{s(j)}\circ\varphi_{i_j}|H=\eta|H\}
$$
is meagre. Define $\theta: \hat A^k\to\hat H$ by
$$
\theta(\chi_1,\ldots,\chi_k)=\prod_{j=1}^{n}\chi_j^{s(j)}\circ\varphi_{i_j}|H.
$$
Then $\theta$ is a continuous homomorphism. In fact $\theta$ is onto
since if $\chi\in\hat H$ then let $\tilde\chi\in\hat A$ be a
character extending $\chi^{s(1)}\circ\varphi_{i_1}^{-1}$ to $A$,
then
$$
\theta(\tilde\chi,1,\ldots, 1)=\chi.
$$
Thus by the open mapping Theorem $\theta$ is open, and so by Pettis'
Theorem if the set
$$
\ker(\theta)=\{(\chi_1,\ldots,\chi_n)\in\hat
A^n:\theta(\chi_1,\ldots,\chi_n)=1\}
$$
is non-meagre then $\{1\}$ is open in $\hat H$, and so $\hat H$ must
be discrete and $H$ must be finite. But this contradicts that
$[A:H]<\infty$. So $\ker(\theta)$ is meagre, and so since
$$
(\eta_1,\ldots,\eta_n)\ker(\theta)=\{(\chi_1,\ldots,\chi_n)\in\hat
A^n:\prod_{j=1}^{n}\chi_j^{s(j)}\circ\varphi_{i_j}|H=\eta|H\}
$$
where $\theta(\eta_1,\ldots,\eta_n)=\eta|H$, it follows that this
set is meagre. \qed (Claim 2)

\bigskip

Applying the Mycielski-Kuratowski Theorem, it follows from Claim 1
and 2 that there is a perfect set $C\subseteq \hat A$ such that for
all $i$ such that $H_i$ is infinite and $j$ such that $H_j'$ is
infinite and all $\chi,\eta\in C$ distinct that $(\chi,\eta)\in
R_{H_i}$ and $(\chi,\eta)\in R_{H_j'}$, and for all
$\chi_1,\ldots,\chi_{k+1}\in C$ distinct we have
$(\chi_1,\ldots,\chi_{k+1})\in Q^{(i_1,\ldots, i_n)}_s$. Then this
$C$ is a perfect sufficiently independent set.
\end{proof}

\bigskip

\begin{proof}[Proof of Theorem 1]
By Theorem 5.1 in \cite{kecsof}, we have that $\approx^{P_c(2^\N)}$
is a generically turbulent equivalence relation. Thus
$E_0<_B\approx_c^{2^\N}$ and $\approx_c$ is not classifiable by
countable structures.

Let $C\subseteq \hat A$ be a perfect sufficiently independent
subset, which we assume has Haar measure 0. For each $\sigma\in
P_c(C)$ let $\hat\sigma$ be the symmetrization and let
$\rho_{\hat\sigma}$ be defined as in 3.4. It is clear that if
$\sigma,\sigma'\in P_c(C)$ are absolutely equivalent then
$\rho_{\hat\sigma}$ and $\rho_{\hat\sigma'}$ are conjugate. If on
the other hand $\sigma\not\approx\sigma'$, then by Lemma 3.10
$\rho_{\hat\sigma}|A$ and $\rho_{\hat\sigma'}|A$ are not conjugate,
since otherwise $\Phi_{\hat\sigma}\approx\Phi_{\hat\sigma'}$.

If $\rho_{\hat\sigma}$ is not a.e. free, then we may instead
consider the product with the Bernoulli shift action on
$(\R^G,\mu_{\varphi_{reg}})$, where $\mu_{\varphi_{reg}}$ is the
product of standard $N(0,1)$ Gaussian probability measures. Then
this action is free and as an $A$-action it has same spectral
properties as $\rho_{\hat\sigma}$.

Since the assignment $P_c(C)\to P_c(\R^G):\sigma\mapsto\mu_{\sigma}$
is continuous we have that
$$
\approx^{P_c(C)}\leq_B(\simeq^{\wmix_A^*(G,X,\mu)},\simeq_A^{\wmix_A^*(G,X,\mu)}).
$$
Since for $A\leqslant H\leqslant G$ we have that
$\simeq^{\wmix_A^*(G,X,\mu)}\subseteq\simeq_H^{\wmix_A^*(G,X,\mu)}\subseteq\simeq_A^{\wmix_A^*(G,X,\mu)}$
the remaining claims of Theorem 1 now follow from Proposition 2.4
and \cite{kecsof} Theorem 5.1.
\end{proof}

\myheadpar{The not type I case.} The argument above unfortunately
depends completely on the fact that Abelian groups have a ``nice''
representation theory. In case $H$ is not abelian by finite, i.e.
not type I, it is known by a result of Hjorth \cite{hjorth1} that
the infinite dimensional irreducible unitary representations of $H$
cannot be classified by countable structures, in fact unitary
equivalence $\sim_u$ in the space $\Irr_\infty(H)$ of infinite
dimensional irreducible representations has a turbulent part. Thus
we have the following:

\medskip

\thm{Theorem 2.}{Suppose $G$ is a countable group and $H$ is an
infinite subgroup which is not abelian-by-finite. Let $(X,\mu)$ be a
non-atomic standard Borel probability space, and let $\mathscr X$ be
an infinite dimensional separable Hilbert space. Then
$\sim_u^{\Irr_\infty(H, \mathscr X)}\leq_B^\omega
\simeq_H^{\wmix^*(G,X)}$. So in particular $\simeq_H^{\wmix^*(G,X)}$
is not classifiable by countable structures and $E_0<_B
\simeq_H^{\wmix^*(G,X)}$}

\begin{proof}
(The proof follows the idea in \cite{zimmer}, p. 110-111. See also
Appendix E of \cite{kechris2} for a detailed account of this
argument.) Let $\pi$ be an irreducible unitary representation of $H$
on $\mathscr X$. Then we can find a weakly mixing Gaussian action
$\rho_{\pi}$ of $G$ such that $\Ind_H^G(\pi)$ is a summand of the
Koopman representation $\kappa$. By Theorem E.1 in \cite{kechris2},
$\pi\mapsto\rho_{\pi}$ maps equivalent unitary representation $\pi$
to conjugate $G$-actions. On the other hand, since $\pi<\kappa|H$ at
most countably many inequivalent $\pi$ can give rise to
conjugate-on-$H$ actions. Thus $\pi\mapsto\rho_{\pi}$ is a
countable-to-1 reduction of $\sim_u$ to $\simeq_H$, and the Theorem
follows from Proposition 2.4.
\end{proof}

Unfortunately the above actions may fail to be ergodic on
$H\leqslant G$, since in general $\Ind_H^G(\pi)|H$ may not be weakly
mixing. If, however, $H\triangleleft G$ is a normal infinite
subgroup and $(\pi,\mathscr X)$ is a weakly mixing representation of
$H$, then $\Ind_H^G(\pi)|H$ is weakly mixing. More generally, we
will say that a subgroup $H\leqslant G$ is {\it index stable} if for
all $g\in G$ we have that
$$
[H:gHg^{-1}\cap H]<\infty \iff [H:g^{-1}Hg\cap H]<\infty.
$$
Then we have:

\medskip

\thm{Theorem 3.}{Suppose $H<G$ is an infinite index stable subgroup
of $G$ which is not abelian by finite. If $\mathscr X$ is an
infinite dimensional separable Hilbert space then
$$
\sim_u^{\Irr_\infty(H, \mathscr X)}\leq_B^\omega
\simeq_H^{\wmix^*_H(G,X)}
$$
}

\begin{proof}
It suffices to show that if $(\pi,\mathscr X)$ is a weakly mixing
unitary representation of $H$ then $\Ind_H^G(\pi)|H$ is weakly
mixing.

We again identify $\Ind_H^G(\pi)$ with $\pi^\alpha$ for some cocycle
$\alpha:G\times G/H\to H$ as defined in \cite{kechris2}, appendix G.
Let $f:G/H\to\mathscr X$ be a unit vector in $l^2(G/H,\mathscr X)$.
We claim that $\langle\pi^\alpha(h)(f):h\in H\rangle$ is not finite
dimensional.

For this, suppose $f(gH)\neq 0$ for some $g=g_i$. Assume first that
 $[H:g^{-1}Hg\cap H]<\infty$. Since for $h\in gHg^{-1}\cap H$ we have $hg=gh'$ for some $h'\in H$
we have $\alpha(h,gH)=g^{-1}hg$ for all $h\in gHg^{-1}\cap H$. So
for $h\in gHg^{-1}\cap H$ we get
$$
\pi^\alpha(h)(f)(gH)=\pi(\alpha(h^{-1},gH)^{-1})(f(h^{-1}gH))=\pi(g^{-1}hg)(f(gH)).
$$
Let $P_{gH}:l^2(G/H,\mathscr X)\to\mathscr X$ be the projection map
$$
P_{gH}(f)=f(gH).
$$
Then the above shows that
$$
P_{gH}(\langle \pi^\alpha(h)(f):h\in gHg^{-1}\cap
H\rangle)=\langle\pi(h)(P_{gH}(f)):h\in g^{-1}Hg\cap H\rangle.
$$
But the latter cannot be finite dimensional since that would
contradict the Lemma 4.4, since we assumed that $[H:g^{-1}Hg\cap
H]<\infty$.

So we may assume that $[H:g^{-1} Hg\cap H]=\infty$ for all cosets
$gH$ where $f(gH)\neq 0$. Suppose for a contradiction that
$\langle\pi^\alpha(H)(f)\rangle=\langle\pi^\alpha(h)(f):h\in
H\rangle$ is finite dimensional, and let $\xi_1,\ldots, \xi_n$ be an
orthonormal basis for it. If $u\in\langle\pi^\alpha(H)(f)\rangle$ is
a unit vector and
$$
u=a_1\xi_1+\cdots+a_n\xi_n
$$
then
\begin{align*}
\|P_{gH}(u)\|&\leq |a_1|\|P_{gH}(\xi_1)\|+\cdots+|a_n|\|P_{gH}(\xi_n)\|\\
&\leq (|a_1|+\cdots+|a_n|)\max_{1\leq k\leq n}\|P_{gH}(\xi_k)\|\\
&\leq n \max_{1\leq k\leq n}\|P_{gH}(\xi_k)\|.
\end{align*}
Let $r=\max_{i} \|f(g_iH)\|$ and let $g_{i_0}$ be such that
$\|f(g_{i_0}H)\|=r$. We have $r>0$ since $f$ is non-zero. Let
$$
\Delta=\{g_i: (\exists k\leq n)\|P_{g_iH}(\xi_k)\|\geq r/n\}.
$$
Then $\Delta$ is finite. It follows from the above that for
$g_i\notin\Delta$ and unit vectors
$u\in\langle\pi^\alpha(H)(f)\rangle$ that
$$
\|P_{g_iH}(u)\|< r.
$$
Now we must have that $h\mapsto hg_{i_0}H$ ranges over infinitely
many cosets: Indeed, if $hg_{i_0}H=h'g_{i_0}H$ then we have
$hg_{i_0}\tilde h=h'g_{i_0}$ for some $\tilde h\in H$ and so
$h^{-1}h'=g_{i_0}\tilde h g_{i_0}^{-1}$. Thus
$$
h^{-1}h'\in g_{i_0} H g_{i_0}^{-1}\cap H.
$$
So if $h,h'\in H$ are not in the same coset of $H/(g_{i_0} H
g_{i_0}^{-1}\cap H)$ then $hg_{i_0}H\neq h'g_{i_0}H$, which proves
that $h\mapsto hg_{i_0}H$ has infinite range since we assumed that
$[H:g_{i_0}^{-1}Hg_{i_0}\cap H]=\infty$, and this implies that
$[H:g_{i_0}Hg_{i_0}^{-1}\cap H]=\infty$.

It follows that we can find $h\in H$ such that $hg_{i_0}H\neq g_iH$
for all $g_i\in\Delta$. But now for this $h$ we have
$$
\|\pi^\alpha(h)(f)(hg_{i_0}H)\|=\|\pi(\alpha(h,hg_{i_0})^{-1})(f(g_{i_0}H))\|=\|f(g_{i_0}H)\|=r
$$
while by the above we also have
$$
\|\pi^\alpha(h)(f)(hg_{i_0}H)\|=\|P_{hg_{i_0}H}(\pi^\alpha(h)(f))\|<r,
$$
since $\pi^\alpha(h)(f)$ is a unit vector in
$\langle\pi^\alpha(H)(f)\rangle$ and $hg_{i_0}\notin\Delta$.
\end{proof}

\mysec{Orbit Equivalence.}

We will now prove an immediate generalization of Corollary 2.6 in
\cite{hjorth2}.

\mythm{Lemma.}{Let $G$ be a countable group with the relative
property (T) over the infinite subgroup $H<G$. Then at most
countably many ergodic-on-$H$ probability measure preserving actions
of $G$ that are not conjugate on $H$ can be orbit equivalent.}

\mypar Before the proof, recall that if $E$ is a measure preserving
countable equivalence relation on a standard Borel probability space
$(X,\mu)$ then the {\it inner group} or {\it full group} of $E$ is
$$
[E]=\{T\in\Aut(X,\mu): (\forall^\mu x) xET(x)\}.
$$
Since $E$ is countable, $[E]$ becomes a Polish group when given the
{\it uniform} topology, which is induced by the metric
$$
d_U(T,S)=\mu\{x\in X: T(x)\neq S(x)\}.
$$
Further, we may define a Borel measure $M$ on $E$ by
$$
M(A)=\int |A_x| d\mu(x).
$$
Where $A_x=\{y: (x,y)\in A\}$. Since $E$ is measure preserving
$$
M(A)=\int |A^y|d\mu(y)
$$
defines the same measure, where $A^y=\{x:(x,y)\in A\}$. The
equivalent integral form is
$$
\int f(x,y)dM(x,y)=\int\sum_{y\in [x]_E} f(x,y)d\mu(x).
$$
For each $T\in[E]$ there is a corresponding unit vector in $L^2(E)$,
$$
\Delta_T(x,y)=\left\{\begin{array}{ll} 1 & \text{ if } T(x)=y\\
0 & \text{ if } T(x)\neq y
\end{array}\right.
$$
and so it is natural to think of the space $L^2(E,M)$ as a
``linearization'' of $[E]$.

\begin{proof}[Proof of Lemma 4.1.]
Let $(Q,\varepsilon)$ be a Kazhdan pair for $(G,H)$ such that any
$(Q,\varepsilon)$-invariant unit vector is within $\frac 1 {10}$ of
an $H$-invariant vector. Let $\sigma\in\erg_H(G,X,\mu)$ and denote
by $\mathcal A_\sigma$ the set of $\tau\in\erg_H(G,X,\mu)$ such that
$E_\tau\subseteq E_\sigma$ (a.e.). We equip $\mathcal A_\sigma$ with
the topology it inherits as a subspace of $[E_\sigma]^G$ when
$[E_\sigma]$ is given the uniform topology. Note that $\mathcal
A_\sigma$ is separable.

We will show that conjugacy on $H$ in $\mathcal A_\sigma$ is an open
equivalence relation. Then it follows by separability that $\mathcal
A_\sigma$ can only meet countably many $\sim_H$ classes, from which
the Lemma follows.

Let $\tau_1,\tau_2\in\mathcal A_\sigma$ and suppose for all $g\in Q$
that
$$
\mu(\{x\in X: \tau_1(g)^{-1}(x)\neq \tau_2(g)^{-1}(x)\})<\frac
{\varepsilon^2}{2}.
$$
We claim that $\tau_1$ and $\tau_2$ are conjugate on $H$. To see
this, define a unitary representation $\pi^{\tau_1,\tau_2}$ on
$L^2(E_\sigma, M)$ by
$$
\pi^{\tau_1,\tau_2}(g)(f)(x,y)=f(\tau_1(g)^{-1}(x),\tau_2(g)^{-1}(y)).
$$
Denote by $I$ the identity transformation. Then clearly
$\pi^{\tau_1,\tau_2}(g)(\Delta_I)=\Delta_{\tau_2(g)\tau_1(g)^{-1}}$,
and so for $g\in Q$ we have from our assumption that
\begin{align*}
&\|\pi^{\tau_1,\tau_2}(g)(\Delta_I)-\Delta_I\|^2_{L^2(E_\sigma)}=\|\Delta_{\tau_2(g)\tau_1(g)^{-1}}-\Delta\|_{L^2(E_\sigma)}^2\\
&=\int \sum_{y\in [x]_{E_\sigma}}
|\Delta_{\tau_2(g)\tau_1(g)^{-1}}(x,y)-\Delta_I(x,y)|^2d\mu(x,y)<\varepsilon^2.
\end{align*}
Thus $\Delta_I$ is $(Q,\varepsilon)$-invariant, and so there is a
$H$-invariant vector $f\in L^2(E_\sigma)$ with
$\|f-\Delta_I\|_{L^2(E_\sigma)}\leq \frac 1 {10}$. Since
$\|f-\Delta_I\|_{L^2(E_\sigma)}\leq \frac 1 {10}$ we must have that
$$
\mu(\{x\in X: |1-f(x,x)|\geq\frac 1 2\})<\frac 1 4.
$$
On the other hand, we must also have that
$$
\mu(\{x\in X: (\exists y_1,y_2) y_1\neq y_2\wedge f(x,y_1)> \frac 1
2\wedge f(x,y_2)> \frac 1 2\})<\frac 1 4.
$$
It follows that the set
$$
A=\{x\in X: (\exists! y) f(x,y)>\frac 1 2\}
$$
has $\mu(A)>0$. Since moreover $f$ is $H$-invariant, $A$ is
$\tau_1|H$-invariant. Since $\tau_1$ is ergodic in $H$ we must have
$\mu(A)=1$.

Define $\varphi(x)$ to be the unique $y\in X$ such that
$f(x,y)>\frac 1 2$. Then clearly $\varphi\subseteq E_\sigma$. It is
also clear from the above argument that $\varphi$ is invertible,
since $\tau_2$ is ergodic on $H$. Thus $\varphi$ is measure
preserving.

From the $H$-invariance of $f$ we now get that for $h\in H$,
$$
\frac 1
2<f(x,\varphi(x))=\pi^{\tau_1,\tau_2}(h)(f)(x,\varphi(x))=f(\tau_1(h)^{-1}(x),\tau_2(h)^{-1}(\varphi(x))).
$$
Since $\varphi(x)$ is the unique $y$ where $f(x,y)>1/2$ it follows
that
$$
\varphi(\tau_1(h)^{-1}(x))=\tau_2(h)^{-1}(\varphi(x))
$$
for all $h\in H$. Thus $\varphi$ witnesses that $\tau_1$ and
$\tau_2$ are conjugate on $H$.
\end{proof}

Let $(X,\mu)$ be a standard Borel probability space. We denote by
$\simeq_{OE}$ the equivalence relation of orbit equivalence in
$\act(G,X,\mu)$. Combining Lemma 4.1 and Theorem 1 and 3, we then
have:

\bigskip

\thm{Theorem 4.}{Let $G$ be a countable group with the relative
property (T) over an infinite subgroup $H$, and let $(X,\mu)$ be a
standard Borel probability space. If either

\(i\) $H$ is Abelian \emph{or}

\(ii\) $H$ is not Abelian-by-finite and is index stable in $G$,

\noindent then $G$ has continuum many orbit inequivalent weakly
mixing a.e. free actions on $(X,\mu)$. In fact, $E_0<_B
\simeq_{OE}^{\wmix^*(G,X,\mu)}$ and $\simeq_{OE}^{\wmix^*(G,X,\mu)}$
is not classifiable by countable structures.}

\begin{proof}
({\it i}) By Theorem 1 we have
$\approx^{P_c(2^\N)}\leq_B(\simeq^{\wmix^*_{A}(G,X,\mu)},\simeq_A^{\wmix^*_{A}(G,X,\mu)})$.
Let $\theta:P_c(2^\N)\to \wmix^*_{A}(G,X,\mu)$ witness this. Then
$\theta$ is also homomorphism of $P_c(2^\N)$ into
$\simeq_{OE}^{\wmix^*(G,X,\mu)}$. Since if $\mu,\mu'\in P_c(2^\N)$
are absolutely inequivalent measures then $\theta(\mu)$ and
$\theta(\mu')$ are not conjugate on $A$, it follows by Lemma 4.1
that at most countably many absolutely inequivalent measures in
$P_c(2^\N)$ can give rise to orbit equivalent actions. Thus $\theta$
is Borel countable-to-1 reduction of absolute equivalence in
$P_c(2^\N)$ to $\simeq_{OE}^{\wmix^*(G,X,\mu)}$. Theorem 4 (i) now
follows from Proposition 2.4 and Theorem 5.1 in \cite{kecsof}.

({\it ii}) Let $\mathscr X$ be an infinite dimensional Hilbert
space. It is clear that the reduction $\pi\mapsto\rho_\pi$ given in
Theorem 2 is a homomorphism of $\sim_u$ in $\Irr_\infty(H,\mathscr
X)$ to $\simeq$ and therefore also to $\simeq_{OE}$. If $H$ is index
stable then by Theorem 3 the action $\rho_{\pi}$ are weakly mixing
on $H$. By Lemma 4.1 if follows that only countably many
inequivalent $\pi$ can give rise to $\rho_\pi$ that are conjugate on
$H$. Therefore the homomorphism $\pi\mapsto\rho_\pi$ is
countable-to-1 reduction of $\sim_u^{\Irr_\infty(H)}$ to
$\simeq_{OE}^{\wmix^*(G,X,\mu)}$.
\end{proof}

\bigskip

\begin{small}
{\sc\noindent Department of Mathematics, University of Toronto\\
40 St. George Street, Room 6092, Toronto, Ontario, Canada\\
{\it E-mail}: {\tt asger@math.utoronto.ca}}
\end{small}


\begin{thebibliography}{3}
\begin{small}

\bibitem{beke} H. Becker, A. Kechris,
{\it The descriptive set theory of Polish group actions}, London
Mathematical Society lecture notes, vol. 232 (1996), Cambridge
University Press.

\bibitem{behava} B. Bekka, P. de la Harpe, A. Valette, {\it Kazhdan's Property
(T)}, to appear, (2003).

\bibitem{bogachev} V.I. Bogachev, {\it Gaussian Measures},
Mathematical Surveys and Monographs, America Mathematical Society
(1998).

\bibitem{cfw} A. Connes, J. Feldman, B. Weiss,
{\it An amenable equivalence relation is generated by a single
transformation}, Ergodic Theory and Dynamical Systems, vol. 1(1981),
pp. 431--450.

\bibitem{connesweiss}
A. Connes, B. Weiss, {\it Property ${\rm T}$ and asymptotically
invariant sequences}, Israel Journal of Mathematics, vol. 37 (1980),
pp. 209--210.

\bibitem{dye} H.A. Dye,
{\it On groups of measure preserving transformation. I}, American
Journal of Mathematics, vol. 81(1959) pp. 119--159.

\bibitem{dye2} H.A. Dye,
{\it On groups of measure preserving transformation. II}, American
Journal of Mathematics, vol. 85(1963) pp. 551--576.

\bibitem{foremanweiss} M. Foreman, B. Weiss, {\it An anti-classification theorem
for ergodic measure preserving transformations}. J. Eur. Math. Soc.
6 (2004), no. 3, pp. 277--292

\bibitem{glasner} Glasner, E., {\it Ergodic Theory via Joinings}, Mathematical Surveys and Monographs
no. 101, American Mathematical Society (2003).

\bibitem{gw} E. Glasner, B. Weiss
{\it Kazhdan's property T and the geometry of the collection of
invariant measures}, GAFA, Geom. Funct. anal. 7, (1997), pp.
917--935.

\bibitem{halmos} P. Halmos,
{\it Lectures on ergodic theory}, Chelsea Publishing Co., New York
(1956).

\bibitem{hjorth1} G. Hjorth,
{\it When it's bad it's worse}, (note), (1997).

\bibitem{hjorth2} G. Hjorth,
{\it A converse to Dye's theorem}, Trans. Amer. Math. Soc., 357
(2005) no. 8, pp. 3083--3103

\bibitem{hjorth3} G. Hjorth, {\it Classification and Orbit Equivalence Relations}, Mathematical Surveys and Monographs 75,
American Mathematical Society (2000)

\bibitem{kechris} A. Kechris,
{\it Classical descriptive set theory}, Springer Verlag, New York,
(1995).

\bibitem{kechris2} A. Kechris,
{\it Global aspects of ergodic group actions and equivalence
relations}, preprint, Caltech, (2006).

\bibitem{kecsof} A. Kechris, N.E. Sofronidis,
{\it A strong ergodicity property of unitary and self-adjoint
operators}, Ergododic Theory and Dynamamical Systems 21, pp.
1459-–1479, (2001).

\bibitem{ornsteinweiss} D. Ornstein, B. Weiss,
{\it Ergodic theory of amenable group actions}, Bulletin of the
American Mathematical Society, vol. 2(1980), pp. 161--164.

\bibitem{popa} Popa, S. {\it Some computations of 1-cohomology groups and
construction of non-orbit-equivalent actions}, J. Inst. Math.
Jussieu 5 (2006), no. 2, pp. 309--332

\bibitem{sinai} I.P. Cornfield, S.V. Fomin, Ya. G. Sinai,
{\it Ergodic Theory}, Springer-Verlag, New York, (1982).

\bibitem{schmidt} K. Schmidt,
{\it Asymptotically invariant sequences and an action of SL$(2, \Z)$
on the sphere}, Israel Journal of Mathematics, vol. 37 (1980), pp.
193--208.

\bibitem{atphd} A. T\"ornquist, {\it The Borel complexity of orbit equivalence}, Ph.D. Thesis, UCLA (2005).

\bibitem{zimmer} R. Zimmer,
{\it Ergodic theory and semisimple Lie groups}, Birkh\"auser,
(1984).

\end{small}

\end{thebibliography}
\end{document}